\documentclass[conference]{IEEEtran}
\IEEEoverridecommandlockouts
\usepackage{amsmath,amssymb,amsfonts,mathrsfs,mathtools,bm, bbm, dsfont,mathrsfs,amsthm,blkarray}
\usepackage[mathscr]{eucal}
\usepackage[dvipdfm]{graphicx}
\usepackage[dvipsnames]{xcolor}
\usepackage{makecell}
\usepackage{float}
\usepackage[font=small]{caption}
\usepackage{epsfig,epsf,psfrag,latexsym, graphics}
\usepackage[ruled,vlined]{algorithm2e}
\usepackage{booktabs}
\usepackage{multirow,multicol}
\usepackage[breaklinks,colorlinks, linkcolor=MidnightBlue, anchorcolor=MidnightBlue, citecolor=MidnightBlue, urlcolor=MidnightBlue]{hyperref}
\usepackage[hyphenbreaks]{breakurl}
\usepackage{cite}
\usepackage{xurl}
\usepackage{threeparttable}
\usepackage{textcomp}

\def\beq{\begin{equation}}
\def\eeq{\end{equation}}
\def\bea{\begin{eqnarray}}
\def\eea{\end{eqnarray}}
\def\ba{\begin{array}}
\def\ea{\end{array}}

\def\bitem{\begin{itemize}}
\def\eitem{\end{itemize}}
\def\ben{\begin{enumerate}}
\def\een{\end{enumerate}}

\def\ie{{\it i.e.,\ \/}}

\newcommand{\beqa}{\begin{eqnarray}}
\newcommand{\eeqa}{\end{eqnarray}}
\newcommand{\beqan}{\begin{eqnarray*}}
\newcommand{\eeqan}{\end{eqnarray*}}

\renewcommand{\[}{\left[}

\newcounter{l1}
\newcounter{l2}
\newcounter{l3}
\newcommand{\bdotlist}{\begin{list}{$\bullet$}{}}
\newcommand{\bboxlist}{\begin{list}{$\Box$}{}}
\newcommand{\bbboxlist}{\begin{list}{\raisebox{.005in}{{\tiny
$\blacksquare$ \ \ }}}{}}
\newcommand{\bdashlist}{\begin{list}{$-$}{} }
\newcommand{\blist}{\begin{list}{}{} }
\newcommand{\barablist}{\begin{list}{\arabic{l1}}{\usecounter{l1}}}
\newcommand{\balphlist}{\begin{list}{(\alph{l2})}{\usecounter{l2}}}
\newcommand{\bAlphlist}{\begin{list}{\Alph{l2}.}{\usecounter{l2}}}
\newcommand{\bdiamlist}{\begin{list}{$\diamond$}{}}
\newcommand{\bromalist}{\begin{list}{(\roman{l3})}{\usecounter{l3}}}

\def\ie{{\it i.e.,\ \/}}
\def\BibTeX{{\rm B\kern-.05em{\sc i\kern-.025em b}\kern-.08em
    T\kern-.1667em\lower.7ex\hbox{E}\kern-.125emX}}

\begin{document}

\title{Energy Management for Renewable-Colocated\\[0.1em] Artificial Intelligence Data Centers
}
\author{\large Siying Li, Lang Tong, and Timothy D. Mount
\thanks{\scriptsize
Siying Li and Lang Tong (\{sl2843, lt35\}@cornell.edu) are with the School of Electrical and Computer Engineering, Cornell University, Ithaca NY, USA. Timothy D. Mount (tdm2@cornell.edu) is with the Dyson School of Applied Economics and Management, Cornell University, Ithaca NY, USA.}
\thanks{\scriptsize This work was supported in part by the National Science Foundation under Awards 2218110 and 2419622, and the Power Systems Engineering Research Center (PSERC).}
}


\maketitle

\begin{abstract}
We develop an energy management system (EMS) for artificial intelligence (AI) data centers with colocated renewable generation. Under a cost-minimizing framework, the EMS of renewable-colocated data center (RCDC) co-optimizes AI workload scheduling, on-site renewable utilization, and electricity market participation. Within both wholesale and retail market participation models, the economic benefit of the RCDC operation is maximized. Empirical evaluations using real-world traces of electricity prices, data center power consumption, and renewable generation demonstrate significant electricity cost reduction from renewable and AI data center colocations.
\end{abstract}


\section{Introduction}
Data centers are among the largest consumers of electricity. While the overall electricity demand in the US increased by 2.3\% over the past decade, the electricity demand from data centers rose by 300\%, now accounting for approximately 4.4\% of total US electricity consumption. The advent of the artificial intelligence (AI) revolution is accelerating the growth of electricity demand. According to Lawrence Berkeley National Laboratory (LBNL) \cite{doe2024datacenters}, electricity demand from US data centers could double or even triple by 2028, increasing data center demand to 6.7\% to 12\% of total US electricity consumption.

The escalating electricity demand from data centers poses
significant challenges for the power grid. Accommodating
this surge requires substantial investments in grid expansion
and generation capacity. However, such expansion typically involves long lead times, with some utilities reporting wait times of up to 7 to 10 years \cite{Duke2025,wecc_load}. From the perspective of the grid operation, colocating data centers with on-site renewable energy sources offers a promising approach to mitigating operational challenges. By enabling the direct use of locally generated renewable power, renewable colocation can reduce the net demand from data centers, alleviate network congestion, and advance decarbonization objectives.

How does renewable colocation benefit data centers? On-site generation of renewable reduces the cost of power consumption and carbon footprint, of course. However, will the cost reduction and environmental benefits offset the investment costs and the complexity of system integration for on-site renewable generation? To this end, it is essential to characterize the benefits of renewable colocation in terms of cost savings, arising from both reduced grid electricity usage and revenues from electricity market participation that offset operational costs, while ensuring that the AI data center prioritizes serving AI tasks.

This paper focuses on the optimal energy management of a renewable-colocated data center (RCDC) by modeling the role of RCDC as a {\em flexible prosumer} participating in a wholesale or retail market while prioritizing the computational needs of its AI customers. The flexibility of AI data centers stems from the characteristics of two types of AI tasks: delay-tolerant (and thus deferrable) AI model-training tasks and delay-intolerant (non-deferrable) AI inferencing tasks, with the former accounting for 30-40\% of annual energy consumption, according to the Electric Power Research Institute (EPRI) \cite{epri2024ai_energy}. Leveraging the inherent flexibility of AI tasks, we develop a scheduling strategy that optimally allocates deferrable AI tasks to specific periods based on real-time renewable
generation and electricity prices, thereby maximizing the economic benefits of RCDCs.

\subsection{Related Work}
Colocation of large load has gained much interest recently \cite{ferc2025colocation,konidena2024colocation}. Major applications include data center colocation \cite{ORO&etal:15,Cao&etal:22} and renewable-colocated manufacturing \cite{Li&etal:25arxiv}. Here, we focus on the energy management of renewable-colocated data centers. A particularly relevant line of work is to align data center workloads with renewable energy availability. For example, Goiri et al. \cite{GOIRI2015520} proposed a scheduling framework in which a data center was powered by solar, with the electrical grid serving as a backup. The proposed energy management system (EMS) predicts short-term solar energy availability and dynamically allocates workloads to maximize the utilization of green energy while ensuring job deadlines are met. Similarly, Google has implemented a workload allocation strategy that shifts computational tasks to times of high renewable generation, thereby advancing toward the goal of operating on carbon-free energy \cite{google2020}.

Exploiting task flexibility in data center workload scheduling has been extensively studied with a primary focus on improving energy efficiency while maintaining quality of service. In response to the growing volume of tasks generated by the digital economy, Yuan et al. \cite{YuanEtal:21} proposed a multi-objective optimization framework to determine task distribution across multiple internet service providers and the service rates of individual data centers. Their approach aims to maximize data center profits while minimizing the average task loss probability across all applications. Other studies have proposed task scheduling methods that minimize the makespan, energy consumption, execution overhead, and the number of active racks \cite{SHARMA2020100373} or jointly optimize energy consumption and execution time \cite{JUAREZ2018257}.

Except for \cite{Li&etal:25arxiv}, which focuses on the renewable-colocated manufacturing of green hydrogen, existing work has overlooked the revenue-enhancing potential of enabling RCDCs to participate in electricity markets as prosumers. While the high-level approach presented here aligns with that of \cite{Li&etal:25arxiv}, there are significant differences between data center operations and green hydrogen production.

\subsection{Contributions}
The results presented here represent the first attempt to capture the full economic benefits of renewable-colocated data centers, featuring bidirectional participation in both wholesale and retail electricity markets. Our goal is to characterize the cost-minimization potential of data center-renewable colocation under a framework that captures the power exchange dynamics between the RCDC and the grid, incorporating the operational characteristics of data center workloads.

The proposed approach has three novel features. First, we model AI data center as a prosumer in its interaction with the wholesale and retail electricity markets, which account for additional data center revenue not considered in existing work.

Second, we present a workload scheduling strategy that leverages the temporal characteristics of AI data center workloads. The EMS incorporates real-time on-site renewable generation and electricity price information, allowing deferrable tasks to be strategically delayed to reduce net costs through optimized task processing and renewable energy sales.

Finally, we conduct numerical studies using real-world datasets, including renewable generation profiles, data center power traces, and electricity prices. Preliminary results show that the proposed optimal EMS can reduce monthly electricity costs by 79.48\% for an RCDC participating in the wholesale market and by 64.53\% in the retail market. These savings exceed the amortized monthly cost of renewable investment, thereby justifying colocation. More comprehensive results are presented in \cite{LiTongMount:arxiv}.
\begin{table}[htbp]
\caption{Variables and System Parameters}
\vspace{-0.1cm}
\begin{center}
\begin{tabular}{|l|l|}
\hline
\multicolumn{2}{|l|}{\textbf{System Parameters}} \\
\hline
$\boldsymbol{\alpha}, \boldsymbol{\beta}$ & Power consumption coefficients of the data center\\
$Q_{\mbox {\tiny R}}$, $Q_{\mbox {\tiny D}}$ & Renewable/data center nameplate capacity \\
\hline
\multicolumn{2}{|l|}{\textbf{Exogenous Variables}} \\
\hline
$\eta_t$ & Capacity factor of renewable generation at time $t$ \\
$\lambda$ & Demand charge \\
$\pi^+_t, \pi^-_t$ & Retail rates for purchasing/selling electricity at time $t$\\
$\pi^{\mbox {\tiny LMP}}_t$ & Wholesale market locational marginal price at time $t$ \\
$W^{\mbox {\tiny DF}}$ & Total deferrable tasks during the scheduling horizon\\
$W_t^{\mbox {\tiny ND}}$ & Non-deferrable data center tasks at time $t$\\
\hline
\multicolumn{2}{|l|}{\textbf{Decision Variables}} \\
\hline
$P_t^{\mbox {\tiny D}}$ & Power input of the data center at time $t$\\
$P_t^{\mbox {\tiny IM}}$, $P_t^{\mbox {\tiny EX}}$ & Power imported from/exported to the grid at time $t$\\
$W_t^{\mbox {\tiny DF}}$ & Deferrable data center tasks processed at time $t$\\
\hline
\end{tabular}
\label{tab:parameters}
\end{center}
\vspace{-0.5cm}
\end{table} 
\subsection{Notations}
The notations used in this paper follow standard conventions. We use $x$ to represent a scalar and $\boldsymbol{x}$ for a vector. We define $[x]:=\{1,\cdots,x\}$. For future reference, key variables and system parameters are summarized in Table~\ref{tab:parameters}.

\section{Renewable-Colocated Data Center}\label{sec:Setting}
With on-site renewable generation, a grid-connected RCDC can draw power from and export excess power to the grid. This section presents the analytical model for the RCDC operation and its interaction with the grid.

\subsection{RCDC Operation Model}\label{sec:RCDC-Model}
\paragraph{Data Center Workload Execution} 
RCDC's workload processing is defined by intervals, indexed by $t=1,\cdots, T$ and of length $\Delta T$\footnote{Without loss of generality, we assume $\Delta T=1$.}.  Within interval $t$, we model RCDC's operation by its ``processing" function $F$ that maps data center's power input $P_t^{\mbox {\tiny D}}$ in kW to computing rate in giga floating-point operations per second (GFLOPS). Function $F$ is nonlinear in general.  Here we assume a piecewise linear approximation in the form:
\begin{equation}\label{eq:computing-rate}
    W_t =F(P_t^{\mbox {\tiny D}})\Delta T =\left(\sum_{k=1}^{K} \left( \alpha_k P_t^{\mbox {\tiny D}}+\beta_k \right) \mathds{1}_{\{P_t^{\mbox {\tiny D}}\in \mathcal{P}_k^{\mbox {\tiny D}}\}}\right)\Delta T,
\end{equation}
where $W_t$ denotes the amount of task processed at time $t$. $K$ denotes the number of segments in the piecewise linear function. Each segment is characterized by a pair of parameters $(\alpha_k,\beta_k)$, with $\{(\alpha_k,\beta_k)\}_{k \in [K]}$ specifying the slope and intercept. The corresponding segment domains are denoted by $\mathcal{P}_k^{\text{\tiny D}}$. The indictor function $\mathds{1}_{\{P_t^{\mbox {\tiny D}}\in \mathcal{P}_k^{\mbox {\tiny D}}\}}$ equals one if the power input falls within segment $k$, and zero otherwise. The piecewise linear formulation, adopted in prior work such as \cite{radovanovic2021powermodeling}, provides a tractable and effective representation of the data center's efficiency profile. In practice, \eqref{eq:computing-rate} can be derived from empirical workload and power consumption data.

\paragraph{RCDC Market Participation}
We analyze a setting where the RCDC participates in either the real-time wholesale electricity market or the retail electricity market. In the wholesale market, electricity imports and exports are settled at the locational marginal price (LMP) $\pi^{\mbox {\tiny LMP}}_t$ (\$/kWh). 

Retail market pricing is regulated under a specific tariff for the RCDC as a commercial customer. A standard tariff applied by most distribution utilities is the net energy metering (NEM) tariff, which prices customer imports and exports differently. In particular, energy drawn from the grid costs the RCDC $\pi_t^+$ (\$/kWh) at $t$, while energy exported to the grid yields a payment of $\pi_t^-$ (\$/kWh). To recover infrastructure-related costs, utilities commonly impose demand charges at $\lambda$ (\$/kW) for RCDCs as large loads. RCDC may also be subject to other fixed (connection) charges, which are ignored in our formulation since they do not affect RCDC scheduling decisions. 

At each time interval $t\in [T]$ over the scheduling horizon $T$, the RCDC makes real-time operational decisions regarding its grid interaction. Specifically, it determines the quantity of power to import from the grid, $P_t^{\mbox {\tiny IM}}$, or to export as excess renewable energy, $P_t^{\mbox {\tiny EX}}$, with both traded at their respective electricity prices. The power exchange between the RCDC and the grid reflects the net balance between the data center's power input $P_t^{\mbox {\tiny D}}$ and the available renewable generation $\eta_tQ_{\mbox {\tiny R}}$.

\vspace{-0.1cm}
\subsection{RCDC Workload Scheduling}
Compared to conventional data centers, AI data centers offer significantly greater flexibility in workload scheduling. Their operations typically involve a combination of non-deferrable and deferrable tasks. Non-deferrable tasks, such as inference requests, require immediate processing to meet user demands. In contrast, tasks like training neural networks for large language models and other machine learning applications are deferrable and can be scheduled based on resource availability and cost considerations \cite{Duke2025}.

We categorize data center AI tasks into these two types. Non-deferrable tasks are executed upon arrival, while deferrable tasks may be scheduled flexibly within their specified deadlines. The cost of processing deferrable tasks varies over time, depending on the availability of renewable generation and the electricity prices. This cost variability creates an opportunity to reduce operational expenses by shifting deferrable tasks to lower-cost intervals, provided that the overall processing requirements are met.

In time interval $t$, the EMS first allocates power to meet the demand of non-deferrable tasks $W_t^{\mbox {\tiny ND}}$. Given the available renewable generation and the electricity price, the RCDC determines the additional power required for deferrable workloads. The total power input $P_t^{\mbox {\tiny D}}$ is the sum of the power allocated to non-deferrable and deferrable tasks.

\vspace{-0.1cm}
\section{Cost-Minimizing Operation}\label{sec:profit-maximization}
The RCDC generates revenue through task processing and by exporting excess renewable power to the grid, while incurring costs from purchasing grid power when on-site generation is insufficient. We assume that the data center does not drop any non-deferrable tasks and that the deadline of deferrable tasks is strict\footnote{For simplicity, we consider a single set of deferrable tasks. The model can be readily extended to accommodate multiple sets with varying deadlines.}. Accordingly, for the purpose of EMS optimization, we exclude the revenue from AI customers, as it is independent of EMS decisions provided that customer requests are satisfied, and focus solely on the operational net cost associated with electricity market interactions.

Specifically, the EMS objective function $J_{\mbox {\tiny T}}(\boldsymbol{P}^{\mbox {\tiny IM}}, \boldsymbol{P}^{\mbox {\tiny EX}};\boldsymbol{\pi})$ is defined as the RCDC's net electricity cost, determined by electricity prices, grid purchases, and surplus renewable exports.
\begin{equation}
    J_{\mbox {\tiny T}}(\boldsymbol{P}^{\mbox {\tiny IM}}, \boldsymbol{P}^{\mbox {\tiny EX}};\boldsymbol{\pi})=J_{\mbox {\tiny T}}^{\mbox {\tiny IM}}(\boldsymbol{P}^{\mbox {\tiny IM}};\boldsymbol{\pi})-J_{\mbox {\tiny T}}^{\mbox {\tiny EX}}(\boldsymbol{P}^{\mbox {\tiny EX}};\boldsymbol{\pi}),
\end{equation}
where $\boldsymbol{P}^{\mbox {\tiny IM}}=\{P_t^{\mbox {\tiny IM}}\}_{t=1}^T$ denotes the power imported from the grid, and $\boldsymbol{P}^{\mbox {\tiny EX}}=\{P_t^{\mbox {\tiny EX}}\}_{t=1}^T$ represents the power exported to the grid. The vector $\boldsymbol{\pi}$ represents the pricing parameters of the wholesale or retail markets, as defined below.

Prices in the wholesale market are determined by the real-time LMPs, \ie$\boldsymbol{\pi}^{\mbox {\tiny WS}}=\{\pi^{\mbox {\tiny LMP}}_t\}_{t=1}^T$. In the retail market, as described in Sec. \ref{sec:RCDC-Model}, the pricing structure includes electricity import and export charges, as well as a demand charge. 
Thus, the retail market pricing is characterized by $\boldsymbol{\pi}^{\mbox {\tiny RT}}=\left( \{ \pi_t^+ \}_{t=1}^T,\ \{ \pi_t^- \}_{t=1}^T,\ \lambda\right)$.

The cost-minimization problem is formulated as follows:
\begin{subequations}\label{eq:profit-maximization}
\begin{align}
& \underset{\{\boldsymbol{P}^{\mbox {\tiny EX}},\boldsymbol{P}^{\mbox {\tiny IM}},\boldsymbol{P}^{\mbox {\tiny D}}, \boldsymbol{W}^{\mbox {\tiny DF}}\}}{\rm minimize} && J_{\mbox {\tiny T}}(\boldsymbol{P}^{\mbox {\tiny IM}},\boldsymbol{P}^{\mbox {\tiny EX}};\boldsymbol{\pi}) \\ \label{eq:total-deferrable}
& \mbox{subject to} && \mathbf{1}^\top \boldsymbol{W}^{\mbox {\tiny DF}}\geq W^{\mbox {\tiny DF}}, \\ \label{eq:power-input}
&&&  \boldsymbol{W}^{\mbox {\tiny DF}}+\boldsymbol{W}^{\mbox {\tiny ND}} =\boldsymbol{F}\left( \boldsymbol{P}^{\mbox {\tiny D}} \right),\\
\label{eq:task-bound}
&&& \boldsymbol{W}^{\mbox {\tiny DF}}\succeq \boldsymbol{0},\\ \label{eq:energy-balance}
&&& \underline{\boldsymbol{P}}^{\mbox {\tiny NET}}\preceq \boldsymbol{P}^{\mbox {\tiny D}}+\boldsymbol{P}^{\mbox {\tiny EX}}-\boldsymbol{P}^{\mbox {\tiny IM}}\preceq \overline{\boldsymbol{P}}^{\mbox {\tiny NET}},\\ \label{eq:power-bound1}
&&& \underline{\boldsymbol{P}}^{\mbox {\tiny D}}\preceq \boldsymbol{P}^{\mbox {\tiny D}}\preceq \overline{\boldsymbol{P}}^{\mbox {\tiny D}},\\
&&& \underline{\boldsymbol{P}}^{\mbox {\tiny EX}}\preceq \boldsymbol{P}^{\mbox {\tiny EX}}\preceq \overline{\boldsymbol{P}}^{\mbox {\tiny EX}},\\ \label{eq:power-bound2}
&&& \underline{\boldsymbol{P}}^{\mbox {\tiny IM}}\preceq \boldsymbol{P}^{\mbox {\tiny IM}}\preceq \overline{\boldsymbol{P}}^{\mbox {\tiny IM}},\\
\label{eq:non-simultaneous}
&&& \boldsymbol{P}^{\mbox {\tiny EX}}\odot\boldsymbol{P}^{\mbox {\tiny IM}}=\boldsymbol{0}.
\end{align}
\end{subequations}

\begin{table*}[htbp]
\centering
\begin{threeparttable}
\caption{Monthly operational outcomes under different configurations and markets}
\begin{tabular}{llccc}
\toprule
\textbf{Market} & \textbf{Quantity} & \textbf{No Colocation} & \textbf{Colocation} & \textbf{Optimal Colocation} \\
\midrule

\multirow{4}{*}{Wholesale} 
  & Electricity Imported (MWh) & $53,221.88$ & $18,980.22$ & $20,283.18$ \\
  & Renewable Exported (MWh)   & --  & $8,897.73$ & $10,200.69$ \\
  & Self Consumption of Renewables (MWh)   & --  & $34,241.66$ & $32,938.70$ \\
  & Electricity Cost (Percentage Reduction\tnote{a} )  & \$$3.82\times 10^6$ & \$$1.19\times 10^6$ (68.8\%) & \$$0.78\times 10^6$ (79.5\%)\\
  & Investment-Adjusted Cost Reduction  & --  & \$$0.17\times 10^6$ & \$$0.58\times 10^6$ \\
\midrule

\multirow{5}{*}{Retail} 
  & Electricity Imported (MWh) & $53,221.88$ & $18,980.22$ & $20,841.47$ \\
  & Renewable Exported (MWh)   & --  & $8,897.73$ & $10,758.97$ \\
  & Self Consumption of Renewables (MWh)   & --  & $34,241.66$ & $32,380.42$\\
  & Peak Demand (kW) & $77,000.00$ & $72,911.64$ & $68,293.31$ \\
  & Electricity Cost (Percentage Reduction\tnote{a} )  & \$$4.88\times 10^6$ & \$$2.04\times 10^6$ (58.2\%) & \$$1.73\times 10^6$ (64.5\%)\\
  & Investment-Adjusted Cost Reduction  & --  & \$$0.38\times 10^6$ & \$$0.69\times 10^6$ \\
\bottomrule
\end{tabular}\label{tab:monthly-total}
\begin{tablenotes}
\footnotesize
\item[a] Percentage of cost reduction compared to the no colocation configuration.
\item[] \vspace{-0.6cm}
\end{tablenotes}
\end{threeparttable}
\end{table*}
In the optimization, constraint \eqref{eq:total-deferrable} ensures that all deferrable tasks are completed within the scheduling horizon. Constraints \eqref{eq:power-input}-\eqref{eq:task-bound} specify the workload execution at each time interval. 
Constraint \eqref{eq:energy-balance} enforces the power balance between renewable generation, grid interaction, and data center consumption. Constraints \eqref{eq:power-bound1}-\eqref{eq:power-bound2} impose limits on power input to the data center, as well as on the amounts of power that can be exported to or imported from the grid. Constraint \eqref{eq:non-simultaneous} prevents simultaneous export and import of grid power. 

This real-time decision-making process can be implemented using a model predictive control (MPC) framework, driven by forecasts of renewable generation and electricity prices. Due to the nonconvexity of the optimization, solving it efficiently is nontrivial; an efficient implementation can be found in \cite{LiTongMount:arxiv}.


\section{Numerical Study}\label{sec:simulation}
We modeled a data center with a nameplate capacity of $Q_{\mbox {\tiny D}}=100$ MW, colocated with a wind farm of capacity $Q_{\mbox {\tiny R}}=150$ MW\footnote{The percentage-based results in this section remain applicable under proportional scaling of both data center and renewable capacities, assuming renewable generation scales linearly with its installed capacity.}. The renewable generation profile and electricity price data were derived from publicly available sources in New York State \cite{NYISO:Data, NationalGrid2025ServiceRates}. In particular, the wholesale prices $\pi^{\mbox {\tiny LMP}}_t$ were based on real-time LMPs, while the retail energy charges $\pi^+_t$ and $\pi^-_t$ were fixed hourly rates tied to day-ahead market prices. The demand charge was set at $\lambda=12.39$\$/kW for a 15-minute peak window.
The workload profile was based on power traces from Google's cluster management software and systems \cite{powerdata2019}, and workload processing was modeled using a two-segment piecewise linear function. Operational decisions for the RCDC were made at 15-minute intervals over a 24-hour horizon and evaluated over a one-month period.

We assessed the operation and economic performance of the following three configurations in both wholesale and retail electricity markets:
\begin{itemize}
    \item No colocation: Grid-only operation without renewable colocation, representing conventional data center operations.
    \item Colocation: Integration of colocated renewables with bidirectional market participation, but without workload scheduling.
    \item Optimal colocation: Full integration of renewable generation, market participation, and workload scheduling, with joint optimization of energy procurement, sales, and workload execution based on real-time electricity prices and renewable availability, as formulated in \eqref{eq:profit-maximization}.
\end{itemize}

\subsection{Colocation Cost Benefits}
Table \ref{tab:monthly-total} summarizes RCDC's monthly operational outcomes under both wholesale and retail settings, assuming 40\% deferrable workload. Compared to the no colocation configuration, we observed 68.8\% cost reduction in the wholesale market and 58.2\% in the retail market under the colocation configuration, while the optimal colocation configuration achieved reductions of 79.5\% and 64.5\%, respectively\footnote{The presented electricity cost is the net cost, which accounts for revenue from selling renewables to the grid.}. According to the wind energy cost report published by the National Renewable Energy Laboratory (NREL), the capital expenditure for utility-scale, land-based wind projects is \$1,968/kW, with operational expenditures of \$43/kW/year in 2023 \cite{NREL91775}. Assuming a 30-year lifespan and a corresponding monthly mortgage interest rate of 0.564\%\footnote{This interest rate reflects typical US commercial bank lending rates \cite{fred_prime_rate}.}, the amortized monthly cost for a 150 MW land-based wind installation is approximately \$2.46 million. The cost reduction from renewable colocation could offset this amount, supporting the economic viability of renewable-colocated data centers.

The benefit of colocation arose from reduced grid electricity purchases, while optimal colocation further enhanced savings by scheduling deferrable tasks and strategically participating in electricity markets. Compared to the simple RCDC, the optimal configuration exhibited greater market participation, buying grid electricity during low-price periods and selling renewables at higher prices. In the retail market, it additionally mitigated peak demand, thereby lowering the demand charge.

\subsection{Impact of Flexibility on Cost Reduction}
\begin{figure}[htbp]
\centerline{\includegraphics[width=0.7\linewidth]{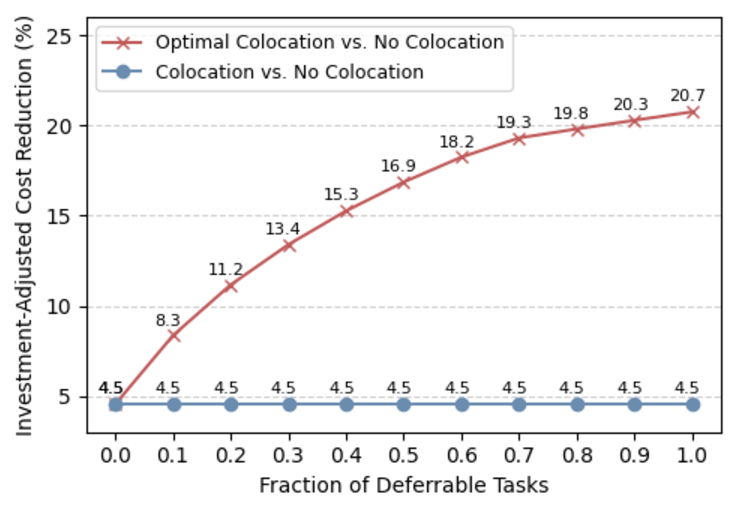}}
\vspace{-0.1cm}
\caption{\scriptsize Investment-adjusted cost reduction in the wholesale market under varying fractions of deferrable tasks.}
\label{fig:impact-fraction}
\end{figure}

\begin{figure}[htbp]
\centerline{\includegraphics[width=0.7\linewidth]{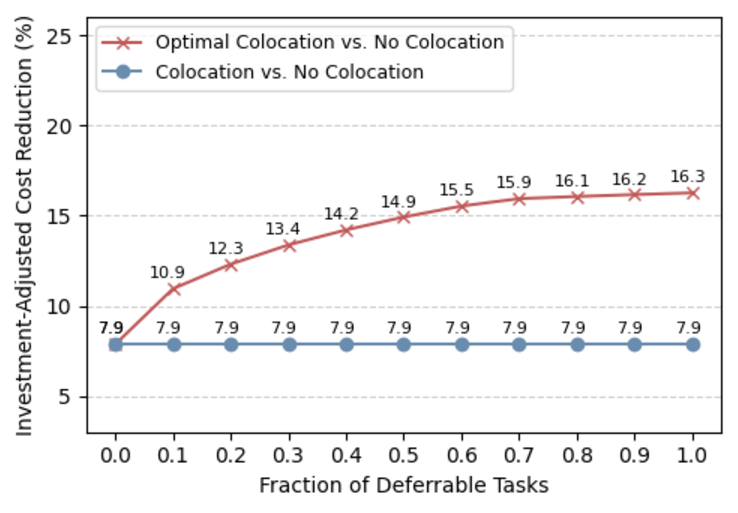}}
\vspace{-0.1cm}
\caption{\scriptsize Investment-adjusted cost reduction in the retail market under varying fractions of deferrable tasks.}
\label{fig:impact-fraction-sc2}
\vspace{-0.1cm}
\end{figure}

Figs. \ref{fig:impact-fraction}-\ref{fig:impact-fraction-sc2} indicate that RCDC's optimal EMS achieved larger cost reductions as the share of deferrable tasks increased, with diminishing marginal benefits at higher fractions. This trend was more pronounced in the wholesale market, whereas cost reductions in the retail market were less sensitive to the deferrable task fraction.
\vspace{-0.2cm}
\vspace{-0.5cm}

In both markets, a higher proportion of deferrable tasks enabled the data center to shift more workload to periods with lower electricity prices, enhancing savings under the optimal colocation configuration. By contrast, in the colocation configuration, cost reduction remained constant regardless of the deferrable task fraction, as workload scheduling was not exploited and economic gains arose solely from renewable energy utilization.

\subsection{Impact of Renewable Generation on Cost Reduction}
\begin{figure}[htbp]
\centerline{\includegraphics[width=0.7\linewidth]{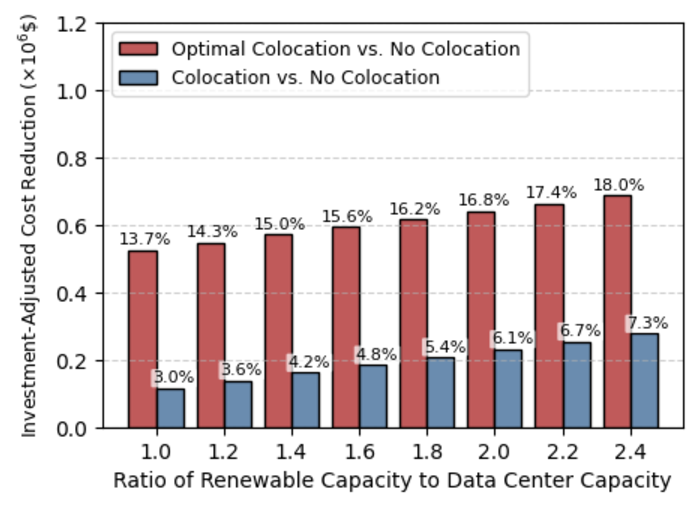}}
\vspace{-0.1cm}
\caption{\scriptsize Investment-adjusted cost reduction in the wholesale market under varying renewable-to-data center capacity ratios.}
\label{fig:ratio-sc1}
\end{figure}

\begin{figure}[htbp]
\centerline{\includegraphics[width=0.7\linewidth]{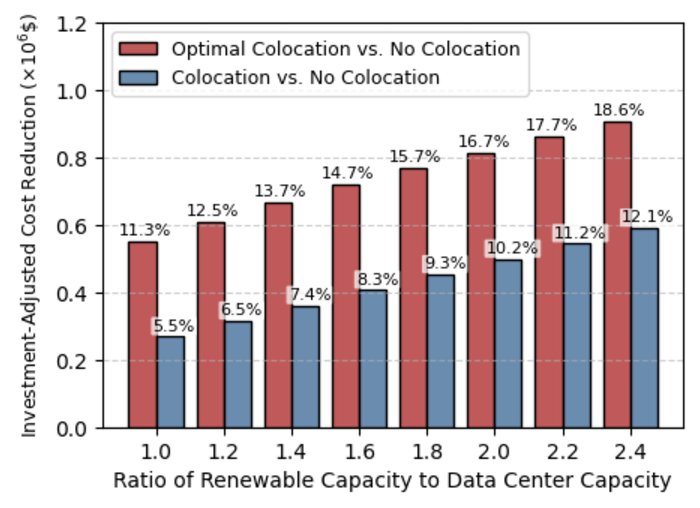}}
\vspace{-0.1cm}
\caption{\scriptsize Investment-adjusted cost reduction in the retail market under varying renewable-to-data center capacity ratios.}
\label{fig:ratio-sc2}
\vspace{-0.1cm}
\end{figure}

Figs. \ref{fig:ratio-sc1}-\ref{fig:ratio-sc2} demonstrate that renewable colocation consistently yielded positive investment-adjusted cost reductions across all renewable-to-data center capacity ratios, compared to the no colocation configuration. Moreover, the cost savings increased with the ratio, indicating that larger colocated renewable capacity strengthened the economic advantage. This effect resulted from higher renewable availability more effectively offset grid purchases. Although larger renewable deployments incurred greater capital investment, the resulting cost reductions proved sufficient to ensure net economic benefits.
\vspace{-0.2cm}
\vspace{-0.6cm}

\section{Conclusion}
Colocating data centers with on-site renewable generation reduces their net electricity demand, alleviates pressure from large load growth on the grid, and contributes to decarbonization. The RCDC model also enables leveraging the temporal flexibility of AI workloads to better align with renewable availability and market dynamics.

We introduce an optimization framework that models the interaction between workload execution and grid power exchange. By strategically scheduling deferrable tasks based on real-time renewable generation and electricity prices, the RCDC can minimize its net electricity costs. As electricity costs become a major component of data center expenses, this approach offers a promising solution for improving cost-efficiency.


{
\bibliographystyle{IEEEtran}
\bibliography{BIB}
}

\end{document}